\documentclass[12pt]{article}
\usepackage{amsmath}
\usepackage{amsthm}
\usepackage{amsfonts}
\usepackage{amssymb}
\usepackage{float}
\usepackage[dvips]{graphicx}
\theoremstyle{plain}
\newtheorem{theorem}{Theorem}[section]
\newtheorem{proposition}[theorem]{Proposition}

\newtheorem{lemma}[theorem]{Lemma}
\newtheorem{fact}[theorem]{Fact}
\newtheorem{question}{Question}[section]
\newtheorem*{questiondash}{Question~\ref{G:compact}$'$}
\theoremstyle{definition}
\newtheorem{definition}[theorem]{Definition}
\newtheorem{example}[theorem]{Example}
\theoremstyle{remark}

\newtheorem{rem}{Remark}

\theoremstyle{plain}
\newtheorem{theoremalph}{Theorem}

\newcounter{thmalph}
\renewcommand{\thethmalph}{\Alph{thmalph}}
\newenvironment{theoremAlph}{\refstepcounter{thmalph}
\medbreak
\noindent
\textbf{Theorem~\thethmalph} \itshape}
{\medbreak}
\makeatletter
\renewcommand{\@makecaption}[2]{%
  \vskip\abovecaptionskip
  \sbox\@tempboxa{#1 #2}%
  \ifdim \wd\@tempboxa >\hsize
    #1: #2\par
  \else
    \global \@minipagefalse
    \hb@xt@\hsize{\hfil\box\@tempboxa\hfil}%
  \fi
  \vskip\belowcaptionskip}
\makeatother
\numberwithin{equation}{section}

\newcommand{\set}[2]{\{{#1}:{#2}\}}

\newcommand{\sesubset}{\rotatebox[origin=B]{-45}{$\subset$}}
\newcommand{\nesubset}{\rotatebox[origin=B]{45}{$\subset$}}
\newcommand{\swsupset}{\rotatebox[origin=B]{45}{$\supset$}}
\newcommand{\nwsupset}{\rotatebox[origin=B]{-45}{$\supset$}}
\begin{document}

\title{A generalized Cartan decomposition for the double coset space
       \large $(U(n_1) \times U(n_2) \times U(n_3))
       \backslash U(n) /(U(p)\times U(q))$}

\author{Toshiyuki Kobayashi
\footnote
{Partly supported
 by Grant-in-Aid for Exploratory Research
 16654014, 
 Japan Society of the Promotion of Science}
\\
Research Institute for Mathematical Sciences,
\\
Kyoto University
}

\date{}

\maketitle

\makeatletter
\renewcommand{\thefootnote}{}
\renewcommand{\@makefntext}[1]{#1}
\makeatother
\footnote{%
\textit{$2000$MSC}: primary  22E4; secondary 32A37,
43A85, 11F67, 53C50, 53D20
\\
\textit{e-mail address}: \texttt{toshi@kurims.kyoto-u.ac.jp}}

\makeatletter
\renewcommand{\l@section}{\@dottedtocline{1}{.8em}{2em}}
\makeatother

\begin{abstract}
Motivated by recent developments
 on visible actions
 on complex manifolds,
 we raise a question
 whether or not the multiplication of three subgroups
 $L$, $G'$ and $H$ surjects a Lie group $G$ 
 in the setting that
 $G/H$ carries a complex structure
 and contains $G'/G' \cap H$
 as a totally real submanifold.  
 
 Particularly important cases are when
 $G/L$ and $G/H$ are generalized flag varieties,
 and we classify pairs of Levi subgroups $(L, H)$
 such that $L G' H = G$, 
 or equivalently,
 the real generalized flag variety $G'/H \cap G'$
 meets every $L$-orbit 
 on the complex generalized flag variety $G/H$
in the setting that $(G, G') = (U(n), O(n))$.  
For such pairs $(L, H)$, we introduce a
 \textit{herringbone stitch} method to find
 a generalized Cartan decomposition
 for the double coset space $L \backslash G/H$,
 for which there has been no general theory in the non-symmetric  case.
Our geometric results provides
 a unified proof
 of various multiplicity-free theorems in representation theory
 of general linear groups.   
\end{abstract}

\tableofcontents

\section{Introduction and statement of main results}
\label{sec:un32}

\noindent
\textbf{1.1.}\enspace
  Our object of study is the double coset space $L \backslash G / H$,
where $L \subset G \supset H$ are a triple of reductive Lie groups.  

In the \lq symmetric case\rq\
(namely, both $(G, L)$ and $(G, H)$ are symmetric pairs),
the theory of the Cartan decomposition
$G=LBH$
 or its variants  gives an explicit description of
 the double coset decomposition $L \backslash G /H$
(e.g., \cite{xfjann, xhoogen, xmatsu, xmatsuki, xmatsucpt, xmaos}).
However,
 in the general case where one of the pairs $(G,L)$ 
and $(G,H)$ is non-symmetric,
 there is no known  structure theory on
the double cosets $L \backslash G /H$
even for a compact Lie group.

Motivated by the recent works of `visible actions'
 (Definition \ref{def:visible})
 on complex manifolds \cite{xkrims, xkvisible} and
multiplicity-free representations 
(e.g.\ the classification of multiplicity-free tensor product
representations of $GL(n)$, see \cite{xkleiden,xstemgl}), 
we have come to realize
 the importance of understanding the double cosets
 $L \backslash G /H$
  in the non-symmetric case such as 
\begin{equation}\label{eqn:1.1}
(G,L,H) = (U(n), U(n_1) \times U(n_2) \times \cdots \times U(n_k),
U(m_1) \times \cdots \times U(m_l)),
\end{equation}
where $n = n_1 + \cdots + n_k = m_1 + \cdots + m_l$.

In this article,
 we  initiate the study of the double cosets
 $L \backslash G /H$
in the \lq{non-symmetric and visible case}\rq\
 by taking \eqref{eqn:1.1}
 as a test case,
and develop new techniques in finding an explicit decomposition for 
$L \backslash G /H$. 

For this,
first we  single out triples that give rise to
 visible actions.
Theorem~\ref{thm:A} gives a 
classification of the triples $(L,G,H)$ such that $G$ has the
decomposition $G=LG'H$ where $G'=O(n)$, 
or equivalently, any $L$-orbit on the \textbf{complex} generalized
flag variety $G/H \simeq \mathcal{B}_{m_1,\dots,m_l}(\mathbb{C}^n)$
(see \eqref{eqn:Bmn}) intersects with its \textbf{real} form 
$\mathcal{B}_{m_1,\dots,m_l}(\mathbb{R}^n)$.
The proof uses an idea of invariant theory arising from quivers.
The classification includes some interesting non-symmetric cases such
as $k=3$, $l=2$
and $\min(n_1+1,n_2+1,n_3+1,m_1,m_2)=2$.
Then, the $L$-action on $G/H$ (and likewise,
the $H$-action on $G/L$) becomes visible in the sense of 
\cite{xkrims} (see Definition~\ref{def:visible}).

Second, we confine ourselves to these triples, 
and prove an analog of the Cartan decomposition
$$G = LBH$$ 
by finding an explicit subset $B$ in $G'$
such that generic points of $B$ form a $J$-transversal totally real slice
 (Definition \ref{def:visible})
 for the $L$-action on $G/H$ of minimal dimension
(see Theorem~\ref{thm:B}).
The novelty of our method in the non-symmetric case is an idea of
`herringbone stitch' (see Section~\ref{subsec:LBH1}).

\medbreak
\noindent
\textbf{1.2.}\enspace
To explain the perspectives of our generalization of the Cartan
decomposition, 
let us recall briefly classic results on
 the double cosets $L \backslash G / H$
in the symmetric case.
A prototype is a theorem due to H. Weyl: 
let $K$ be a connected compact Lie group,
and $T$ a maximal toral subgroup.
Then,
\begin{equation}\label{eqn:maxT}
\mbox{\textit{%
any element of\/ $K$ is conjugate to an element of\/ $T$.
}}
\end{equation}
We set
 $G: = K \times K$,
$A := \set{(t, t^{-1})}{t \in T}$
and identify $K$ with the subgroup
$\operatorname{diag}(K) :=  
 \set{(k, k)}{k \in K}$
of $G$.
Then,
  the statement \eqref{eqn:maxT}
is equivalent to the double coset decomposition: 
\begin{equation}\label{eqn:1.2}
    G = K A K.
\end{equation}

In the above case,
$(G,K) = (K \times K, \operatorname{diag}(K))$
forms a compact symmetric pair.
More generally,
the decomposition \eqref{eqn:1.2} still holds for a Riemannian
symmetric pair $(G,K)$ by taking
$A \simeq \mathbb{T}^k$
($G/K$: compact type) or
$A \simeq \mathbb{R}^k$
($G/K$: non-compact type)
where $k = \operatorname{rank} G/K$.
 Such a decomposition is known as the Cartan decomposition
for symmetric spaces.  
  
 A further generalization of the Cartan decomposition
has been developed over the decades
 under the hypothesis that both $(G,L)$ and $(G,H)$ are
symmetric pairs.
   For example,
Hoogenboom \cite{xhoogen} gave an analog of the Cartan decomposition for
$L\backslash G/H = (U(l)\times U(m))\backslash U(n)/(U(p)\times U(q))$
$(l+m=p+q=n)$ by finding a toral subgroup 
$\mathbb{T}^k$ as its representatives,
where
$k = \min(l,m,p,q)$.
This result is  generalized by 
Matsuki \cite{xmatsucpt}, showing that there exists a toral subgroup $B$
in $G$ such that 
\begin{equation}\label{eqn:1.4}
G = LBH
\end{equation}
 if $G$ is compact
(see Fact~\ref{fact:symmCartan}).
Analogous decomposition also holds in the case where
$G$ is a non-compact reductive Lie group
and $L$ is its maximal compact subgroup,
by taking  a non-compact abelian subgroup
$B$ of dimension 
$\operatorname{rank}_{\mathbb{R}} G/H$
(see Flensted-Jensen \cite{xfjann}).

\medbreak
\noindent
\textbf{1.3.}\enspace
Before explaining a new direction of study in the non-symmetric case, 
we pin down some remarkable
 aspects on the Cartan decomposition in the symmetric case
from algebraic, geometric, and analytic viewpoints.

 Algebraically,
   finding nice representatives of the double coset
    is relevant to the reduction theory,
or the theory of normal forms.
For example,
     the Cartan decomposition \eqref{eqn:1.2} for
     $(G, K) = (GL(n, \mathbb{R}), O(n))$ corresponds to the diagonalization
      of symmetric matrices by orthogonal transformations.
The case $G=G' \times G', L = H = \operatorname{diag}(G')$
with $G' = GL(n,\mathbb{C})$ is equivalent to the theory of Jordan
   normal  forms.

 Geometrically,
 \eqref{eqn:1.4}
   means that every $L$-orbit on the 
 (pseudo-)Riemannian 
   symmetric space $G/H$ meets the flat totally geodesic submanifold
$B/B \cap H$.
The decomposition \eqref{eqn:1.4} is also used in the construction of
 a $G$-equivariant compactification of the 
 symmetric space $G/H$
(see \cite{xborelji} for a survey on various compactifications).
    
 Analytically, the Cartan decomposition 
is particularly important in the analysis of
    asymptotic behavior of global solutions to $G$-invariant
   differential equations on the symmetric space $G/H$
(e.g.\ \cite{xfjann,xhoogen,xkinv}).

\medbreak
\noindent
\textbf{1.4.}\enspace
Now, we consider the non-symmetric case 
$L \subset G \supset H$.
Unlike the symmetric case,
 we cannot expect the existence of an
 abelian subgroup $B$
   such that $LB H$
contains an interior point of $G$ in general,
as is easily observed by the argument of dimensions
(e.g.\ \cite[Introduction]{xkinv}). 
Instead, we raise here the following question:
 
\begin{question}\label{G:compact}
 Does there exist a `nice' subgroup $G'$ such that $ L G' H$
 contains an open subset of $G$?
\end{question}

 For  a compact $G$,
one may strengthen Question~\ref{G:compact}
as follows:
\begin{questiondash}
 Does there exist a `nice' subgroup $G'$ such that $G = L G' H$?
\end{questiondash}
  
The decomposition $G=LG'H$ means that the double coset space
$L \backslash G / H$ can be controlled by a subgroup $G'$.
  
What is a `nice' subgroup $G'$?
In contrast to the previous case that $G/H$ carries a $G$-invariant
Riemannian structure and that the abelian subgroup
$G' = B$ (see \eqref{eqn:1.4}) gives a flat totally geodesic
submanifold of $G/H$,
we are interested in the case that 
$G/H$ carries a $G$-invariant complex structure and that
the subgroup $G'$ gives a 
\textit{totally real} submanifold $G'/G' \cap H$ of 
 $G/H$. 
In the latter case,
the $L$-action on $G/H$ is said to be \textit{previsible}
(\cite[Definition~3.1.1]{xkrims})
if $LG'H$ contains an open subset of $G$.

\medbreak
\noindent
\textbf{1.5.}\enspace
Let us state our main results.
Suppose we are
  in the setting  \eqref{eqn:1.1}
and consider Question \ref{G:compact}$'$ for
$$
    G' := O(n).
$$
In this setting, 
we shall give a necessary and sufficient condition for the multiplication
map $L \times G' \times H \to G$ to be surjective.

In order to clarify its geometric meaning,
let $\mathcal{B}_{m_1,\ldots,m_l} (\mathbb{C}^n)$ denote the complex
(generalized) flag variety: 
\begin{align}
\{ (V_1,\dots,V_l) : \{0\} &= V_0 \subset
   V_1 \subset V_2 \subset \dots \subset V_{l-1} \subset 
  V_l = \mathbb{C}^n,
\nonumber
\\
   \dim V_i &= m_1 + \dots + m_i
  \quad (1 \le i \le l) \}.
\label{eqn:Bmn}
\end{align}
Likewise, the real (generalized) flag variety
$\mathcal{B}_{m_1,\dots,m_l} (\mathbb{R}^n)$
is defined and becomes a totally real submanifold of
$\mathcal{B}_{m_1,\dots,m_l} (\mathbb{C}^n)$.

\begin{theoremalph}
\label{thm:A}
Let $k,l \ge 2$ and
$n = n_1 + \dots + n_k = m_1 + \dots + m_l$
be partitions of $n$ by positive integers.
Let
$$
(G,L,H) = (U(n), U(n_1) \times U(n_2) \times \cdots \times
 U(n_k), U(m_1) \times \cdots \times U(m_l)).
$$
We set $N := \min(n_1,\dots,n_k)$
and $M := \min(m_1,\dots,m_l)$.
Then the following five conditions are equivalent:
\begin{enumerate}
    \renewcommand{\theenumi}{\roman{enumi}}
    \renewcommand{\labelenumi}{\upshape\theenumi)}
\item 
$G = L G' H$.
Here, $G' := O(n)$.
\item 
$\mathcal{B}_{m_1,\dots,m_l} (\mathbb{R}^n) \times
 \mathcal{B}_{n_1,\dots,n_k} (\mathbb{R}^n)$
meets every $G$-orbit on 
$\mathcal{B}_{m_1,\dots,m_l} (\mathbb{C}^n) \times
 \mathcal{B}_{n_1,\dots,n_k} (\mathbb{C}^n)$
by the diagonal action.
\item[\upshape ii)$'$]
$\mathcal{B}_{m_1,\dots,m_l} (\mathbb{R}^n)$
meets every
$L$-orbit on 
$\mathcal{B}_{m_1,\dots,m_l} (\mathbb{C}^n)$.
\item[\upshape ii)$''$]
$\mathcal{B}_{n_1,\dots,n_k} (\mathbb{R}^n)$
meets every $H$-orbit on
$\mathcal{B}_{n_1,\dots,n_k} (\mathbb{C}^n)$.
\item 
One of the following conditions holds:

\begin{tabular}{lllll}
\hbox to 6mm{\upshape 0)\hfill}
&   $k = 2$, && $l = 2$,
\\
\hbox to 6mm{\upshape I)\hfill}
& $k=3$, &$N=1$, &$l=2$,
\\
\hbox to 6mm{\upshape II)\hfill}
& $k=3$, &$N \ge 2$, & $l=2$, & $M=2$,
\\
\hbox to 6mm{\upshape III)\hfill}
& && $l=2$, & $M=1$,
\\
\hbox to 6mm{\upshape I$'$)\hfill}
& $k=2$, && $l=3$, & $M=1$,
\\
\hbox to 6mm{\upshape II$'$)\hfill}
& $k=2$, & $N=2$, & $l=3$, & $M\ge 2$,
\\
\hbox to 6mm{\upshape III$'$)\hfill}
& $k=2$, & $N=1$.
\end{tabular}
\end{enumerate}
\end{theoremalph}

\begin{rem}
Both $(G,L)$ and $(G,H)$ are symmetric pairs if and only if
$(k,l) = (2,2)$, 
namely, 
$(G,L,H)$ is in Case 0.  
\end{rem}

\begin{rem}
\label{rem:1.5.2}
The condition (ii) implies
 that the $L$-action on $G/H$
 is previsible 
 (see Definition \ref{def:visible}).
We shall see in Section \ref{sec:visible}
 that this action is (strongly) visible,
 too
(see also \cite[Corollary~17]{xkrims}).  
\end{rem}

\begin{rem}\label{rem:spherical}
A holomorphic action of a complex reductive group on a complex
manifold $D$ is called \textit{spherical} if its Borel subgroup has an
open orbit on $D$.
We shall see in Section~\ref{sec:concl} that the condition (ii) 
in Theorem~\ref{thm:A} is equivalent to:
\begin{itemize}
\item[\upshape iv)]
$
\mathcal{B}_{m_1,\dots,m_l} (\mathbb{C}^n)
   \times \mathcal{B}_{n_1,\dots,n_k}(\mathbb{C}^n)
$
is a spherical variety of
$
G_{\mathbb{C}}:= 
   GL(n,\mathbb{C})
$.
\end{itemize}
See Littelmann \cite{xlittelmann}
 for the statement (iv)
in the case $k=l=2$,
namely, Case~0 in (iii).
\end{rem}

\begin{rem}\label{rem:polar}
An isometric action of a compact Lie group $L$ on a Riemannian manifold
is called \textit{polar} if there exists a submanifold that meets
every $L$-orbit orthogonally.
Among Cases 0$\sim$III$'$
 in Theorem~\ref{thm:A}, 
 the $L$-action on $G/H$ is polar if and only if $(G,L,H)$ is in
Case 0 (see \cite{xbigo}).
\end{rem}

\medbreak
\noindent
\textbf{1.6.}\enspace
Suppose one of (therefore, all of) the equivalent 
conditions in Theorem~\ref{thm:A} is satisfied.
As a finer structural result of the double coset decomposition
$L \backslash G / H$,
we shall construct a fairly simple subset $B$
of $G' = O(n)$ such that the multiplication map
$L \times B \times H \to G$ is still surjective,
according to Cases 0 $\sim$ III 
 of Theorem~\ref{thm:A}.
We omit Cases~I$'$ $\sim$ III$'$ below because these are essentially
the same with Cases~I $\sim$ III.
For Case~I below,
we may and do assume $n_1 = 1$ without loss of generalities.

\addtocounter{thmalph}{1}
\begin{theoremAlph}
\label{thm:B}
(generalized Cartan decomposition). 
Let 
$$
(G,L,H)=(U(n), U(n_1) \times \cdots \times U(n_k), U(p) \times U(q)).
$$
Then, there exists 
$B \subset O(n)$
such that $G = LBH$, 
where $B$ is of the following form: 
\begin{equation*}
  B \simeq
  \begin{cases}
    \mathbb{T}^{\min(n_1,n_2,p,q)}          &\mbox{\rm(Case 0),}\\
    \mathbb{T}^{\min(p,q,n_2,n_3)} 
    \cdot \mathbb{T}^{\min(p,q,n_2+1,n_3+1)}&\mbox{\rm(Case I),}\\
    \mathbb{T}^2\cdot\mathbb{T}  
                \cdot\mathbb{T}^2           &\mbox{\rm(Case II),}\\
    \underbrace{\mathbb{T} \cdot \ \cdots \ \cdot \mathbb{T}}_{k-1}
                                            &\mbox{\rm(Case III).}\\
  \end{cases}
\end{equation*}
Here, $\mathbb{T}^a\cdot\mathbb{T}^b$
means  
a subset of the form
$\{xy \in G: x \in \mathbb{T}^a, y \in \mathbb{T}^b \}$
for some
toral subgroups
$\mathbb{T}^a$ and $\mathbb{T}^b$.
\end{theoremAlph}

We note that $B$ is no longer a subgroup of
$G$ in Cases I, II and III.

\medbreak
\noindent
\textbf{1.7.}\enspace
This article is organized as follows.
First, we give a proof of
 Theorem~\ref{thm:B} 
(generalized Cartan decomposition $G=LBH$)
by constructing explicitly 
the subset $B$ of $O(n)$.
This is done in
Theorems \ref{thm:LBH0} (Case 0), 
\ref{thm:LBH1} (Case I), \ref{thm:LBH2} (Case II),
and \ref{thm:LBH3} (Case~III),
respectively by using an idea of herringbone stitch.
This also gives a proof of the implication
(iii) $\Rightarrow$ (i)
in Theorem~\ref{thm:A}.
The remaining implications of Theorem~\ref{thm:A} are proved in
Section~\ref{sec:5}.
An application to representation theory is discussed in
Section~\ref{sec:concl}.

Theorem~\ref{thm:A} was announced and used in
\cite[Theorem~3.1]{xkleiden} and \cite[Theorem~16]{xkrims}
and its proof was postponed until this article.
Theorem~\ref{thm:B} was presented in the Oberwolfach workshop on
``Finite and Infinite Dimensional Complex Geometry and Representation
Theory'', 
organized by A. T. Huckleberry, K.-H. Neeb,
and J. A. Wolf, February 2004.
The author thanks the organizers for a wonderful and stimulating
atmosphere of the workshop.

\section{Symmetric case}\label{subsec:LBH0}

This section reviews a well-known fact on the Cartan decomposition for the
symmetric case.
The results here will be used in the non-symmetric case
(Sections \ref{subsec:LBH1}, \ref{subsec:LBH2}, 
and \ref{subsec:LBH3}) 
as a `stitch' (see Diagram~\ref{dia:3.1}, for example). 
Theorem \ref{thm:LBH0} corresponds to
 Theorem~\ref{thm:B} in Case~0, which is proved here.

First, we recall from 
\cite[Theorem~6.10]{xhoogen},
                 \cite[Theorem~1]{xmatsuki}
the following:

\begin{fact}
\label{fact:symmCartan}
Let $G$ be a connected, compact Lie group with Lie algebra\/
$\mathfrak{g}$. 
Suppose that $\tau$ and $\theta$ are two involutive automorphisms of
$G$, 
and that $L$ and $H$ are open subgroups 
of $G^\tau$ and $G^\theta$,
respectively.
We take a maximal abelian subspace\/ $\mathfrak{b}$ in
$\mathfrak{g}^{-\tau,-\theta} 
 := \{ X \in \mathfrak{g}: \tau X = \theta X = -X \}$,
and write $B$ for the connected abelian subgroup with Lie algebra\/
 $\mathfrak{b}$.
Then
$
G = L B H
$.
\end{fact}

Now, let us consider the setting:
\begin{equation}
\label{eqn:n22}
(G,L,H) = (U(n),U(n_1) \times U(n_2), U(p) \times U(q)),
\end{equation}
where $n=n_1+n_2=p+q$.
Then, both $(G,L)$ and $(G,H)$ are symmetric pairs. 
In fact, if we set
$ I_{p,q} := \operatorname{diag} (1,\dots,1,-1,\dots,-1)$
and define an involution $\theta$ by
$\theta(g) := I_{p,q} g I_{p,q}{}^{-1}$,
then $H = G^\theta$.
Likewise, $L = G^\tau$ if we set
$\tau(g) := I_{n_1,n_2} g I_{n_1,n_2}{}^{-1}$.

We set
\begin{equation}
\label{eqn:lsymm}
l:=\min(n_1, n_2, p, q) 
\end{equation}
and define
an abelian subspace: 
$$
\mathfrak{b} :=\sum_{i=1}^l
 \mathbb R (E_{i, n+1-i} - E_{n+1-i,i}).
$$
Then,
$\mathfrak{b}$ is a maximal abelian subspace in
$\mathfrak{g}^{-\tau,-\theta}$,
and $B:=\exp(\mathfrak b)$ is a toral
   subgroup of $O(n)$.
Now, applying Fact~\ref{fact:symmCartan}, we obtain:

\begin{theorem}
\label{thm:LBH0}
 $G=L B H$.
\end{theorem}

\section{Non-symmetric case 1: $\min(n_1, n_2, n_3)=1$}\label{subsec:LBH1}

In Sections~\ref{subsec:LBH1}, \ref{subsec:LBH2}, and
\ref{subsec:LBH3}, 
we give an explicit decomposition formula for the double coset space 
$L\backslash G/H$ in Cases~I, II, and III of Theorem~\ref{thm:B},
respectively, and complete the proof of Theorem~\ref{thm:B},
and therefore that of the implication (iii) $\Rightarrow$ (ii) of
Theorem~\ref{thm:A}.

The distinguishing feature of these three sections is that we are
dealing with the \textbf{non-symmetric pair} $(G,L)$,
for which there is no known general theory on the double coset decomposition
$L\backslash G/H$
of a compact Lie group $G$.
We shall introduce 
a method of \textit{herringbone stitch}
(see Diagram~\ref{dia:3.1})
consisting of symmetric triples 
$G_{2i+1} \subset G_{2i} \supset H_i$
and
$L_i \subset G_{2i+1} \supset G_{2i+2}$
$(i=0,1,2,\ldots)$ 
such that the iteration of the double coset decomposition of 
$G_{2i+1}\backslash G_{2i}/H_i$
and
$L_i\backslash G_{2i+1}/G_{2i+2}$
keeps on toward
a finer structure of 
$L\backslash G / H = L_0 \backslash G_0 / H_0$.

This section treats the most interesting case for Theorem~\ref{thm:B},
namely, 
\begin{equation}
\label{eqn:n32}
(G,L,H) = (U(n), U(n_1)\times U(n_2)\times U(n_3),
   U(p) \times U(q))
\end{equation}
where
 $\min(n_1, n_2, n_3)=1$.
Theorem~\ref{thm:LBH1} below corresponds to 
 Case~I of Theorem~\ref{thm:B}.

Without loss of generality, we may and do assume $n_1=1$.
Thus, $L = U(1) \times U(n_2) \times U(n_3)$
 ($n_2+n_3=n-1$).
We set
\begin{equation}\label{eqn:rankLH}
l:=\min(2p, 2q,2 n_2+1,2 n_3+1).
\end{equation}
A simple computation shows
 $\left[\frac{l}{2}\right] = \min(p,q,n_2,n_3)$
and $\left[\frac{l+1}{2}\right] = \min(p,q,n_2+1,n_3+1)$.
We define two abelian subspaces: 
\begin{alignat*}{1}
\mathfrak b' &:= \sum_{i=1}^{[\frac{l+1}2]}
                 \mathbb R (E_{i, n+1-i} - E_{n+1-i,i}),
\\
\mathfrak b'' &:= \sum_{i=1}^{[\frac{l}2]}
  \mathbb R (E_{i+1, n+1-i} - E_{n+1-i,i+1}).
\end{alignat*}
Then $B':=\exp \mathfrak b'$ and $B'' := \exp \mathfrak b''$
 are  toral subgroups
of dimension
$[\frac{l+1}{2}]$
and $[\frac{l}{2}]$, respectively.
We define a subset of $O(n)$ by
\begin{equation}
B:= B'' B'.
\end{equation}
For the sake of simplicity,
 we shall write also 
 $B= {\mathbb {T}}^{[\frac l 2]}\cdot {\mathbb {T}}^{[\frac {l+1} 2]}$.  
We note that $B$ is a compact manifold of dimension
 $l=[\frac{l}{2}]+[\frac{l+1}2]$
 because $B$
  is diffeomorphic to the homogeneous space
  $(B' \times B'')/(B' \cap B'')$
 and $B' \cap B''$ is a finite subgroup.

We are ready to describe the double coset decomposition for
$L \backslash G / H$
in the case \eqref{eqn:n32}.

\begin{theorem}[generalized Cartan decomposition]\label{thm:LBH1}
$G=L B H$,
 where $B \simeq {\mathbb {T}}^{\min(p,q,n_2,n_3)}\cdot {\mathbb {T}}^{\min(p,q,n_2+1,n_3+1)}$.
\end{theorem}

\begin{proof}
First, for $i \ge 1$,
we define one dimensional toral subgroups by
\begin{alignat*}{1}
B_i &:= \exp \mathbb R (E_{i, n+1-i} - E_{n+1-i,i}),
\\
C_i &:= \exp \mathbb R (E_{i+1, n+1-i} - E_{n+1-i,i+1}).
\end{alignat*}
Then,
$B_1, B_2, \dotsb,$ and $B_{[\frac{l+1}{2}]}$ commute with each other,
 and we have
$$
B'= B_1 B_2 \cdots B_{[\frac{l+1}{2}]}
=B_{[\frac{l+1}{2}]} \dotsb B_2 B_1.
$$

Likewise,
 $B''= C_1 C_2 \dotsb C_{[\frac{l}2]}$.

For $m \ge 1$, we define an embedding of a
one dimensional torus into $G$ by
$$
  \iota_m : \mathbb T \to G,
            \quad          a \mapsto
          \operatorname{diag}
          (\underbrace{a, \dots, a}_{[\frac{m+1}2]},
           \underbrace{1, \dots, 1}_{n-m},
          \underbrace{a, \dots, a}_{[\frac{m}2]}).
$$
We write $T_m$ for its image,
 and define a subgroup $G_m$ of $G$ by
$$
G_m:= T_m \times U(n-m).
$$
Here,
 we regard $U(n-m)$ as a subgroup of $G$
 by identifying with 
$\{I_{[\frac{m+1}{2}]}\} \times U(n-m) \times \{ I_{[\frac{m}{2}]}\}$.
($I_m$ stands for the unit matrix of degree $m$.)
It is convenient to set $T_0=\{e\}$ and $G_0=G$.
Then,
 we have a decreasing sequence of subgroups:
$$
    G=G_0 \supset G_1 \supset G_2 \supset \dotsb.
$$
The point of our definition of $G_m$ is that
\newline\indent
$G_{2i}$ commutes with all of 
$B_1, \dotsb, B_i, C_1, \dotsb, C_{i-1}$.
\newline\indent
$G_{2i+1}$ commutes with all of
$B_1, \dotsb, B_i, C_1, \dotsb, C_{i}$.

\medskip

Next, for $i \ge 0$, we define the following subgroups:
\begin{alignat*}{2}
H_i &:= H \cap G_{2i}
 &&\simeq T_{2i} \times U(p-i) \times U(q-i),
\\
L_i &:= L \cap G_{2i+1}
 &&\simeq T_{2i+1} \times U(n_2-i) \times U(n_3-i).
\end{alignat*}

The following obvious properties play a crucial role in the inductive
step below. 
\begin{alignat}{3}
&H = H_0 \supset H_1 \supset H_2 \supset \cdots; \quad
&&\text{$H_i$ }&&\text{commutes with $B_1, \dots, B_i$,}
\label{eqn:Hseq}
\\
&L = L_0 \supset L_1 \supset L_2 \supset \cdots; 
&&\text{$L_i$ }&&\text{commutes with $C_1, \dots, C_i$.}
\label{eqn:Lseq}
\end{alignat}

With these preparations, let us proceed the proof of Theorem~\ref{thm:LBH1}
along a \textit{herringbone stitch} consisting of
 triples
 $(G_{2i+1}, G_{2i}, H_i)$
 and  $(L_i, G_{2i+1}, G_{2i+2})$
 ($i=0,1,2,\cdots$):
%
\begin{figure}[H]
$$
   \arraycolsep=2pt
\begin{array}{ccccccccccccccccccccc}
   && L_0 &&&& L_1 &&&& L_2 &&&& L_3 &&&&&&\cdots \\
   &&& \sesubset &&&& \sesubset &&&& \sesubset &&&& \sesubset \\
     &&&& G_1 &&&& G_3 &&&& G_5 &&&& G_7 &&&& \cdots \\
  &&&\swsupset && \nwsupset &&\swsupset && \nwsupset 
              &&\swsupset && \nwsupset &&\swsupset && \nwsupset \\
   && G_0 &&&& G_2 &&&& G_4 &&&& G_6 &&&& G_8 &\enspace&\cdots \\
 &\nesubset &&&& \nesubset &&&&\nesubset &&&& \nesubset &&&& \nesubset \\
H_0 &&&& H_1 &&&& H_2 &&&& H_3 &&&& H_4 &&&& \cdots
\end{array}
$$
\renewcommand{\figurename}{Diagram}
\caption{}
\label{dia:3.1}
\end{figure}

We claim that each triple has the following decomposition formula:
\begin{alignat}{2}
&G_{2i} &&=G_{2i+1} B_{i+1} H_i,
\label{eqn:Geven}
\\
&G_{2i+1} &&= L_i C_{i+1} G_{2i+2}.
\label{eqn:Godd}
\end{alignat}
To see this, we first take 
away the trivial factor from $G_{2i}$ and $G_{2i+1}$, respectively.
Then, the following bijections hold
 for $i = 0, 1, \dots$,
\begin{alignat*}{1}
G_{2i+1} \backslash G_{2i}/ H_i
&\simeq 
(U(1)\times U(n-2i-1)) \backslash U(n-2i)/(U(p-i) \times U(q-i)), 
\\
L_i \backslash G_{2i+1}/ G_{2i+2}
&\simeq
(U(n_2-i) \times U(n_3-i))\backslash
 U(n-2i-1)/(U(n-2i-2) \times U(1)).
\end{alignat*}

Since 
the right-hand side is the double coset space by symmetric subgroups,
 we can apply Theorem~\ref{thm:LBH0}. 
Thus, \eqref{eqn:Geven} and \eqref{eqn:Godd} have been proved.

By using \eqref{eqn:Geven} and \eqref{eqn:Godd} 
 iteratively, together with the commutating properties
\eqref{eqn:Hseq} and \eqref{eqn:Lseq},
 we obtain 
\begin{alignat}{1}
G=G_0 &= G_1 B_1 H_0
\nonumber
\\
      & = (L_0 C_1 G_2) B_1 H
\nonumber
\\
      & = L C_1 (G_3 B_2 H_1) B_1 H
\nonumber
\\
      & = L C_1 G_3 B_2 B_1 H
\nonumber
\\
      &=\dotsb
\nonumber
\\
     &= L C_1 \cdots C_i G_{2i} B_i \cdots B_1 H
\label{eqn:evenfin}
\\
    & =L C_1 \cdots C_i G_{2i+1} B_{i+1} \cdots B_1 H.  
\label{eqn:oddfin}
\end{alignat}
If $\min(p,q) \le \min(n_2,n_3)$,
then this equation terminates with
$G_{2i} = H_i$ when $i$ reaches $\min(p,q)$.
Then, $i = [\frac{l}{2}] = [\frac{l+1}{2}]$, 
and we have $G = LB''B'H_iH = LBH$ from \eqref{eqn:evenfin}.
If $\min(p,q) > \min(n_2,n_3)$,
then this equation terminates with
$G_{2i+1} = L_i$ when $i$ reaches $\min(n_2,n_3)$.
Then, $i = [\frac{l}{2}]$ and
$i+1 = [\frac{l+1}{2}]$,
and we have
$G = L L_i B''B'H = LBH$
from \eqref{eqn:oddfin}.
Hence, we have completed the proof of Theorem.
\end{proof}

\section{Non-symmetric case 2: $\min(n_1,n_2,n_3)\ge 2$}\label{subsec:LBH2}

In this section we study the double coset space
$L\backslash G/H$ in Case~II of Theorem~\ref{thm:B},
that is, the  non-symmetric case: 
$$
(G,L,H) = (U(n), U(n_1) \times U(n_2) \times U(n_3),
   U(p) \times U(q)),
$$
where
$n = n_1 + n_2 + n_3 = p+q$,
$\min(n_1,n_2,n_3) \ge 2$
and
$\min(p,q) = 2$.
Without loss of generality,
 we may and do
assume $p = 2$.

We define
 three abelian subspaces: 
\begin{alignat*}{2}
& \mathfrak a_1 
&&{}:= \mathbb R (E_{1, n} - E_{n,1})
                +\mathbb R (E_{2,n-1} - E_{n-1,2}),
\\
& \mathfrak a' 
&&{}:= \mathbb R (E_{n_1+1,n-1}- E_{n-1,n_1+1}),
\\
& \mathfrak a_2 
&&{}:= \mathbb R (E_{n_1+1, n} - E_{n,n_1+1})
                +\mathbb R (E_{n_1+2,n-1} - E_{n-1,n_1+2}),
\end{alignat*}
and correspondingly three
 toral subgroups by 
     $A_1:=\exp \mathfrak a_1$,  
     $A' := \exp \mathfrak a'$,
 and $A_2:=\exp \mathfrak a_2$.
We then set 
\begin{equation}
B:= A_2 A' A_1.
\end{equation}
Note that $B$ is a five dimensional 
subset of $O(n)$, but is no more a subgroup.

Here is a generalized Cartan decomposition for
$L \backslash G/H$
in the non-symmetric setting 
$\min(n_1,n_2,n_3) \ge 2$
and
$\min(p,q) = 2$:

\begin{theorem}\label{thm:LBH2}
$G=L B H.$

\end{theorem}

\begin{proof}
The proof again uses herringbone `stitch', 
of which each stitch is a special case of
Theorem~\ref{thm:LBH1} or Theorem~\ref{thm:LBH0},
respectively.

Step 1.
We define a subgroup $G_1$ of $G$ by 
$$
   G_1 := U(n_1) \times U(n_2+n_3).
$$
Since both $(G, H)$ and $(G, G_1)$ are symmetric pairs,
 we have
\begin{equation}\label{eqn:3-1}
G= G_1 A_1 H
\end{equation}
 by Theorem~\ref{thm:LBH0}.
We observe from \eqref{eqn:lsymm} that
 $A_1$ is of dimension $\min(2, n-2, n_1, n_2+n_3)=2$.

Step 2.
We define a subgroup $H_1$ of $G$ by
$$H_1:= Z_{H \cap G_1}(A_1)
 \simeq \Delta(\mathbb{T}^2) \times U(n_1-2) \times U(n_2+n_3-2).
$$
Here, 
$\Delta(\mathbb{T}^2)
:=\set{\operatorname{diag}(a,
b,1,\cdots, 1, b,a)}{a, b\in \mathbb T}.
$
Then,
taking away the first factor inclusion $U(n_1)$, we have
$$
  L\backslash G_1/H_1
  \simeq (U(n_2)\times U(n_3))\backslash U(n_2+n_3)
  /(U(n_2+n_3-2)\times U(1) \times U(1)).
$$
Applying Theorem~\ref{thm:LBH1} to the right-hand side,
 we obtain
\begin{equation}\label{eqn:3-2}
  G_1 = L A_2 A' H_1.
\end{equation}
We note that $A_2 A'$
 is of dimension $\min(2n_2, 2n_3, 2(n_2+n_3-2)+1, 3)=3$,
 as explained in Section \ref{subsec:LBH1}.

Combining \eqref{eqn:3-1} and \eqref{eqn:3-2},
 we obtain
\begin{alignat*}{1}
G & = G_1 A_1 H
\\
  & = (L A_2 A' H_1) A_1 H
\\
  & = L A_2 A' A_1 H
\\
  & = L B H
\end{alignat*}
because $H_1$ commutes with $A_1$.
Thus, we have shown Theorem~\ref{thm:LBH2}.
\end{proof}

\section{Non-symmetric case 3: $\min(p,q) = 1$}
\label{subsec:LBH3}
This section treats the double coset space
$L \backslash G / H$ in Case~III of Theorem~\ref{thm:B},
that is, the non-symmetric case:
$$
(G,L,H) = (U(n), U(n_1) \times \dots \times U(n_k),
           U(1) \times U(n-1))
$$
for an arbitrary partition $n = n_1 + \dots + n_k$.
In this case,
although the pair $(G,L)$ is non-symmetric,
the symmetric pair $(G,H)$ gives a very simple homogeneous space, 
namely, $G/H \simeq \mathbb{P}^{n-1} \mathbb{C}$.
Thus, it is much easier to find an explicit decomposition
$G=LBH$ than the previous cases in Sections~\ref{subsec:LBH1} and
\ref{subsec:LBH2}.

For $1 \le i \le k-1$, we set
\begin{align*}
H_i &:= -E_{1,n_1+\dots+n_i+1} + E_{n_1+\dots+n_i+1,1}
\\
B_i &:= \exp (\mathbb{R} H_i)  \quad (\simeq \mathbb{T}),
\\
\intertext{and define a $(k-1)$-dimensional
  subset $B$ in $G' = O(n)$ by}
B   &:= B_1 \cdots B_{k-1}.
\end{align*}
We note that $B$ is not a group if $k \ge 3$.
The subset $B$
 is contained in the subgroup $O(k)$
of $O(n)$ in an obvious sense.

Then, we have

\begin{theorem}
\label{thm:LBH3}
$G = LBH$,
where $B = \underbrace{\mathbb{T} \cdot \ \cdots\ \cdot \mathbb{T}}_{k-1}$
$(\subset O(n))$.
\end{theorem}

\begin{proof}
We shall work on the $L$-action on
 $G/H \simeq \mathbb{P}^{n-1} \mathbb{C}$.
First, we observe that the map
$$
 U(n_i) \times \mathbb{R} \to \mathbb{C}^{n_i},
 \quad
 (h_i,a_i) \mapsto h_i\, {}^t\! (a_i, 0, \ldots, 0)
$$
is surjective for each $i$ $(1 \le i \le k)$.
Therefore, the following map
\begin{align}
& U(n_1) \times \cdots \times U(n_k) \times \mathbb{P}^{k-1} \mathbb{R}
     \to \mathbb{P}^{n-1} \mathbb{C},
\label{eqn:LPRPC}
\\
& ((h_1, \ldots, h_k), [a_1 : \cdots : a_k])
     \mapsto [h_1\, {}^t\!(a_1, 0,\ldots, 0) : \cdots : 
              h_k\, {}^t\!(a_k, 0,\ldots,
     0)]
\nonumber
\end{align}
is also surjective.

Next, 
an elementary matrix computation shows
\begin{align*}
&\exp (\theta_1 H_1) \cdots \exp (\theta_{k-1} H_{k-1})
  \, {}^t \! (1,0,\dots,0)
\\
= {}
&{}^t\!(a_1,\dots,0,  a_2,\dots,0,   \ldots, a_k, 0,\dots,0),
\end{align*}
where
$$
\begin{pmatrix}
   a_1 \\ a_2 \\ a_3 \\ \vdots \\ a_{k-1} \\ a_k 
\end{pmatrix}
=
\begin{pmatrix}
   \cos \theta_1  \cos \theta_2 & \cdots
        & \cos \theta_{k-2}  \cos \theta_{k-1}
 \\
   \sin \theta_1  \cos \theta_2 & \cdots
        & \cos \theta_{k-2}  \cos \theta_{k-1}
 \\
   \phantom{\sin \theta_1} \sin \theta_2 & \cdots
        & \cos \theta_{k-2} \cos \theta_{k-1}
 \\
        &\kern-1em \ddots & \vdots
 \\
        && \sin \theta_{k-2}  \cos \theta_{k-1}
 \\
        && \phantom{\sin \theta_{k-2}} \sin \theta_{k-1}
\end{pmatrix}
 \in \mathbb{R}^k \backslash \{0\}.
$$
Hence,
$B \cdot [1:0:\cdots:0] = \mathbb{P}^{k-1}\mathbb{R}$.

Combining with the surjective map
$L \times \mathbb{P}^{k-1} \mathbb{R} \rightarrow G/H$,
we have shown that
 the multiplication map
$
L \times B \times H \to G
$
is surjective.
\end{proof}

Now, Theorems~\ref{thm:LBH0}, \ref{thm:LBH1},
\ref{thm:LBH2} and \ref{thm:LBH3}
show that we have completed the proof of Theorem~\ref{thm:B}.

\section{Proof of Theorem \ref{thm:A}} \label{sec:5}

This section gives a proof of Theorem~\ref{thm:A}.

Suppose we are in the setting of Theorem~\ref{thm:A}.
The equivalence (i) $\Leftrightarrow$ (ii)
(likewise, (i) $\Leftrightarrow$ (ii)$'$ 
and (i)  $\Leftrightarrow$ (ii)$''$)
is clear from the following natural identifications:
$$
G'/G'\cap H \simeq \mathcal{B}_{m_1,\dots,m_l} (\mathbb{R}^n), 
\quad
G/H \simeq \mathcal{B}_{m_1,\dots,m_l} (\mathbb{C}^n).
$$
Further, Theorem~\ref{thm:B} shows (iii) $\Rightarrow$ (i).

Thus, the rest of this section is devoted to the proof of the implication 
(i) $\Rightarrow$ (iii).

We begin with a question about when
 the multiplication map
\begin{equation}\label{eqn:LGH}
L \times G' \times H \to G,
\quad (l,g',h) \mapsto lg'h
\end{equation}
is surjective,
in the general setting that $G = U(n)$, $G' = O(n)$,
and $H$ is a Levi subgroup $G$.
Our key machinery to find a necessary condition for the surjectivity
of \eqref{eqn:LGH} is
 Lemma~\ref{lem:machinenotreal}. 
Let us explain briefly the ideas of 
our strategy that manages
 the three non-commutative subgroups $L$, $G'$ and $H$:
\begin{alignat*}{2}
&L &&\cdots\  \text{finding $L$-invariants (invariant theory),}
\\
&G' &&\cdots\ \text{using the geometric property ($G'$ gives
 real points in $G/H$),}
\\
&H &&\cdots\ \text{realizing $H$ as the isotropy subgroup of the $G$-action.}
\end{alignat*}

For this, we take
$J \in M(n,\mathbb{R})$ and consider the adjoint orbit:
$$
G / G_J \simeq
\{ \operatorname{Ad} (g) J : g \in G \}
\subset M(n,\mathbb{C}),
$$
where
$\operatorname{Ad} (g) J := gJg^{-1}$ and
$$
G_J := \{ g \in G : gJ = Jg \}.
$$
Later, we shall choose $J$ such that
$H$ is conjugate to $G_J$ by an element of $G'$.
Here, we note that the surjectivity of the map \eqref{eqn:LGH} remains
unchanged if we replace $H$ with $aHa^{-1}$ 
and $L$ with $bHb^{-1}$ $(a,b \in G')$.
Then, the following observation:
$$
G' G_J / G_J
 \simeq \{ \operatorname{Ad} (g) J : g \in G' \}
 \subset M(n,\mathbb{R})
$$
will be used in Lemma~\ref{lem:LgJ}
(`management' of $G' \simeq O(n)$),
while an invariant theory will be used in Lemma~\ref{lem:charpol}
(`management' of $L \simeq U(n_1) \times \dots \times U(n_k)$).

\begin{lemma}\label{lem:LgJ}
Let $J \in M(n,\mathbb{R})$,
$G' = O(n)$,
and $L$ a subgroup of $G = U(n)$.
If there exists $g \in G$ such that
\begin{equation}\label{eqn:LgJ}
\operatorname{Ad}(L) (\operatorname{Ad}(g)J) \cap M(n,\mathbb{R}) = \emptyset,
\end{equation}
then $G \supsetneqq L G' G_J$.
\end{lemma}

\begin{proof}
First we observe that
$$
\operatorname{Ad}(G' G_J) J 
= \operatorname{Ad}(G')J \subset M(n,\mathbb{R}).
$$
Then, the condition \eqref{eqn:LgJ} implies
$\operatorname{Ad}(L g) J \cap \operatorname{Ad}(G' G_J) J = \emptyset$,
whence
$Lg \cap G' G_J = \emptyset$.
Therefore,
 $g \notin L G' G_J$.
\end{proof}

Next, we 
fix a partition $n=n_1+\cdots+n_k$,
and
find a sufficient condition for \eqref{eqn:LgJ} in the
setting:
$$
L = U(n_1) \times \dots \times U(n_k).
$$

For this, we fix $l \ge 2$ and take  a loop
$i_0 \to i_1 \to \cdots \to i_l$ 
consisting of non-negative integers $1 , \ldots, k$
such that
\begin{equation}\label{eqn:string}
 i_0 = i_l, \quad
 i_{a-1} \ne i_a \quad
 (a = 1,2,\dots,l).
\end{equation}
Now, let us introduce a non-linear map:
$$
A_{i_0 \cdots i_l} : M(n,\mathbb{C}) \to M(n_{i_0},\mathbb{C})
$$
as follows:
Let $P \in M(n,\mathbb{C})$, and we
 write
$P$ as  
$(P_{ij})_{1 \le i,j \le k}$
in block matrix form such that the $(i,j)$-block
$
P_{ij} \in M(n_i, n_j; \mathbb{C})
$.
We set
$\widetilde{P}_{ij} \in M(n_i, n_j; \mathbb{C})$ by
$$
\widetilde{P}_{ij} := \begin{cases}
                     P_{ij}    & (i < j), \\
                     P^*_{ji}  & (i > j).
          \end{cases}
$$
Then,
$A_{i_0 \cdots i_l}(P)$ is defined by
\begin{equation}\label{eqn:APstring}
A_{i_0 \cdots i_l} (P) := \widetilde{P}_{i_0 i_1} \widetilde{P}_{i_1 i_2}
   \cdots \widetilde{P}_{i_{l-1} i_l}.
\end{equation}

\begin{lemma}  \label{lem:charpol}
If there exists a loop 
$i_0 \to i_1 \to \cdots \to i_l$ $(= i_0)$ such that
at least one of the coefficients of
the characteristic polynomial\/
$\det(\lambda I_{n_{i_0}} - A_{i_0 \cdots i_l} (P) )$
is not real,
 then
$$
\operatorname{Ad} (L) P \cap M(n, \mathbb{R}) = \emptyset.
$$
\end{lemma}

Later, we shall take $P$ to be $\operatorname{Ad}(g)J$ and 
apply this lemma to the following loops:
\begin{enumerate}
\item[\rm{1)}]
$1 \to 2 \to 3 \to 1$,
\begin{equation}\label{eqn:1231}
A_{1231} = P_{12} P_{23} P^*_{13}
\end{equation}
\item[\rm{2)}]
$1 \to 3 \to 2 \to 4 \to 1$,
\begin{equation}\label{eqn:13241}
A_{13241} = P_{13} P^*_{23} P_{24} P^*_{14}
\end{equation}
\item[\rm{3)}]
$1 \to 4 \to 2 \to 4 \to 3 \to 4 \to 1$,
\begin{equation}\label{eqn:1424341}
A_{1424341} = P_{14} P^*_{24} P_{24} P^*_{34} P_{34} P^*_{14}
\end{equation}
\end{enumerate}

\begin{proof}
For a block diagonal matrix 
$l=\begin{pmatrix}
       l_1 \\  & l_2 \\  && \ddots \\ &&& l_k
     \end{pmatrix}     
\in L$,
the transform
$P \mapsto \operatorname{Ad}(l)P$
induces that of the $(i,j)$-block matrix:
$$
  P_{ij} \mapsto l_i P_{ij} l_j^{-1}
\quad
(1 \le i,j \le k). 
$$
Then 
$P_{ij}^*$
is transformed as
$
P_{ij}^* \mapsto (l_i P_{ij} l_j^{-1})^*
   = l_j P_{ij}^* l_i^{-1}
$. 
Hence,
$\widetilde{P}_{ij}$ is transformed as
$$
\widetilde{P}_{ij} \mapsto l_i \widetilde{P}_{ij} l_j^{-1}
\quad (1 \le i, j \le k),
$$
and then
$A_{i_0\cdots i_l}(P)$ is transformed as
$$
A_{i_0 \cdots i_l} (P) \mapsto
   l_{i_0} A_{i_0 \cdots i_l} (P) l_{i_0}^{-1}.
$$
Therefore,
 the characteristic polynomial of
$A_{i_0\cdots i_l} (P) $
is invariant under the transformation 
$P \mapsto \operatorname{Ad}(l)P$.
In particular, if
$\operatorname{Ad}(L)P \cap M(n, \mathbb{R}) \ne \emptyset$,
then
$\det (\lambda I_{n_{i_0}} - A_{i_0\cdots i_l} (P) )
 \in \mathbb{R} [\lambda]$.
By contraposition, Lemma~\ref{lem:charpol} follows.
\end{proof}

Here is a key machinery to show the implication (i) $\Rightarrow$ (iii) 
in Theorem~\ref{thm:A}:

\begin{lemma}\label{lem:machinenotreal}
Let $n = n_1 + \dots + n_k$ be a partition,
and $L = U(n_1) \times \dots \times U(n_k)$
be the natural subgroup of $G = U(n)$.
Suppose $J$ is of a block diagonal matrix:
\begin{equation} \label{eqn:J}
J :=
 \begin{pmatrix}
      J_1 \\ & J_2 \\  &&  \ddots \\ &&& J_k 
 \end{pmatrix}
 \in  M(n, \mathbb{R}), 
\end{equation}
where $J_i \in M(n_i,\mathbb {R})$
$(1 \le i \le k)$.
If there exist a skew Hermitian matrix
$X \in \mathfrak{u}(n)$ and a loop
$i_0 \to i_1 \to \dots \to i_l$ $(= i_0)$
\textup{(}see \eqref{eqn:string}\textup{)} such that
$$
\det (\lambda I_{n_{i_0}} - A_{i_0\cdots i_l} ([X,J]))
  \notin \mathbb{R} [\lambda],
$$
then
the multiplication map
$L \times G' \times G_J \to G$
is not surjective.
Here, we recall $G' = O(n)$.
\end{lemma}

\begin{proof}
We set
$P(\varepsilon) := \operatorname{Ad} (\exp (\varepsilon X))J$. 
In view of Lemmas~\ref{lem:LgJ} and \ref{lem:charpol},
it is sufficient to show
$$
\det (\lambda I_{n_{i_0}} - A_{i_0 \cdots i_l} (P(\varepsilon)))
 \notin \mathbb{R} [\lambda]
$$
for some $\varepsilon > 0$.
We set
$Q := [X,J]$.
The matrix $P(\varepsilon)$ depends real analytically on $\varepsilon$,
and we have
$$
P(\varepsilon)
 = J + \varepsilon Q + O(\varepsilon^2),
$$
as $\varepsilon$ tends to $0$.
In particular,
the $(i,j)$-block matrix
$P_{ij} (\varepsilon)$
$( \in M(n_i, n_j; \mathbb{C}))$
satisfies
$$
P_{ij} (\varepsilon) = \varepsilon Q_{ij} + O(\varepsilon^2)
\quad (\varepsilon \to 0)
$$
for $i \ne j$.
Then, we have
\begin{align}
    &\det (\lambda I_{n_{i_0}} - A_{i_0 \cdots i_l} (P(\varepsilon)))
\nonumber
\\
={} &\det (\lambda I_{n_{i_0}} - \varepsilon^l \widetilde{Q}_{i_0i_1}
        \cdots \widetilde{Q}_{i_{l-1}i_l} +
       O(\varepsilon^{l+1}))
\nonumber
\\
={} &\det (\lambda I_{n_{i_0}} - \varepsilon^l A_{i_0 \cdots i_l}
        (Q) +
       O(\varepsilon^{l+1}))
\nonumber
\\
={} & \sum_{r=0}^{n_{i_0}} \lambda^{n_{i_0}-r} \varepsilon^{rl}
         h_r (\varepsilon),
\label{eqn:her}
\end{align}
where $h_r (\varepsilon)$ $(0 \le r \le n_{i_0})$
are real analytic functions of $\varepsilon$ such that
$$
\det (\lambda I_{n_{i_0}} - A_{i_0\cdots i_l} (Q))
  = \sum_{r=0}^{n_{i_0}} \lambda^{n_{i_0}-r} h_r (0).
$$
From our assumption,
this polynomial is not of real coefficients, namely,
there exists $r$ such that
$h_r(0) \notin \mathbb{R}$.
It follows from \eqref{eqn:her} that
$\det (\lambda I_{n_{i_0}} - A_{i_0 \cdots i_l} (P(\varepsilon)))
  \notin \mathbb{R} [\lambda]$
for any sufficiently small $\varepsilon > 0$.
Hence, we have shown Lemma.
\end{proof}

For the applications of Lemma~\ref{lem:machinenotreal} below,
we shall take a specific choice of a skew Hermitian matrix
$X \in \mathfrak{u}(n)$.
According to the partition 
$n=n_1+\cdots+n_k$,
we write 
$X = (X_{ij})_{1 \le i,j \le k}$
as a block form.
We note that the $(i,j)$ block of
$Q = [X,J]$ is given by
\begin{equation}
\label{eqn:Pe}
Q_{ij} = X_{ij} J_j - J_i X_{ij} \; .
\end{equation}
We also note that if $J \in M(n,\mathbb{R})$ is a diagonal
matrix whose entries consist of non-negative integers
$1,2,\ldots,l$ such that
$$
\# \{ 1 \le a \le n: J_{aa} = i \} = m_i
\quad (1 \le i \le l),
$$
then $G_J$ is conjugate to 
$U(m_1) \times U(m_2) \times \cdots \times U(m_l)$
by an element of $G' = O(n)$.
Since the surjectivity of \eqref{eqn:LGH} remains unchanged
if we take the conjugation of $H$ by $G'$,
we can apply Lemma \ref{lem:machinenotreal} to study the surjectivity
of \eqref{eqn:LGH} in the setting \eqref{eqn:1.1}.

Now, we  apply Lemma~\ref{lem:machinenotreal}
to show the following four propositions:

\begin{proposition}\label{prop:notreal33}
Let $n = n_1 + n_2 + n_3 = m_1 + m_2 + m_3$
be partitions of $n$ by positive integers.
We define (natural) subgroups $L$ and $H$ of $G = U(n)$ by
\begin{align*}
& L := U(n_1) \times U(n_2) \times U(n_3),
\\
& H := U(m_1) \times U(m_2) \times U(m_3),
\end{align*}
and $G' := O(n)$.
Then, $G \supsetneqq L G' H$.
\end{proposition}

In the following three propositions,
we set
$$
(G,L,H) = (U(n), U(n_1) \times \cdots \times U(n_k), U(p) \times U(q)) 
\quad\text{and}\quad G' = O(n)
$$
where $n = n_1 + n_2 + \dots + n_k = p+q$.

\begin{proposition}  \label{prop:notreal}
$G \supsetneqq L G' H$
if
\begin{equation}  \label{eqn:large}
   \min(p,q) \ge 3,
\quad
k=3 \quad\mbox{and} \quad \min(n_1, n_2, n_3) \ge 2.
\end{equation}
\end{proposition}
\begin{proposition}  \label{prop:notreal2}
$G \supsetneqq L G' H$
if
\begin{equation}\label{eqn:large2}  
   \min(p,q) = 2
\quad\mbox{and}\quad
k\ge 4 .
\end{equation}
\end{proposition}
\begin{proposition}  \label{prop:notreal3}
$G \supsetneqq L G' H$
if
\begin{equation}  \label{eqn:large3}
  \min(p,q) \ge 3
\quad\mbox{and}\quad
k \ge 4 .
\end{equation}
\end{proposition}

\begin{proof}[Proof of Proposition~\ref{prop:notreal33}]
We consider the partition
 $n = n_1 + n_2 +n_3$ and the loop
$1 \to 2 \to 3 \to 1$.
We take
 $J = \operatorname{diag}(j_1,\dots,j_n) \in M(n,\mathbb{R})$ 
to be a diagonal matrix with the following
two properties:
\begin{align*}
& j_{1} = 1, \quad
  j_{n_1} = 2, \quad
  j_{n_1+n_2+1} = 3,
\\
& \# \{ k: j_{k} = i \} = m_i \quad
  (1 \le i \le 3).
\end{align*}
Then, the isotropy subgroup
$G_J$ is conjugate to $H$ by an element of $O(n)$.
We fix $z \in \mathbb{C}$
and define a skew Hermitian matrix
$X = (X_{ij})_{1 \le i, j \le 3} \in \mathfrak{u}(n)$
as follows:
$$
X_{12} = \begin{pmatrix}1\\ &\mbox{\Large $0$}\end{pmatrix}, \;
X_{23} = \begin{pmatrix}1\\ &\mbox{\Large $0$}\end{pmatrix}, \;
X_{13} = \begin{pmatrix}z\\ &\mbox{\Large $0$}\end{pmatrix}.
$$
Then it follows from \eqref{eqn:Pe} that
$Q = [X,J]$ has the following block entries:
$$
Q_{12} = \begin{pmatrix}1\\ &\mbox{\Large $0$}\end{pmatrix}, \;
Q_{23} = \begin{pmatrix}1\\ &\mbox{\Large $0$}\end{pmatrix}, \;
Q_{13}^* = \begin{pmatrix}2\overline{z}\\ &\mbox{\Large $0$}\end{pmatrix},
$$
and therefore we have
$$
A_{1231} (Q)
= Q_{12} Q_{23} Q_{13}^*
= \begin{pmatrix} 2\overline{z} \\ & \mbox{\Large $0$}\end{pmatrix}.
$$

Thus,
$
\det (\lambda I_{n_1} - A_{1231} (Q))
   = \lambda^{n_1} - 2\overline{z} \lambda^{n_1 -1}
$.
This does not have real coefficients
if we take $z \notin \mathbb{R}$.
Hence, Proposition follows from Lemma~\ref{lem:machinenotreal}.
\end{proof}

\begin{proof}[Proof of Proposition \ref{prop:notreal}]
We shall apply Lemma~\ref{lem:machinenotreal} with $k = 3$ and the
loop $1 \to 2 \to 3 \to 1$.
In light of the assumption
 \eqref{eqn:large},
an elementary consideration shows that there exist
positive integers
$p_i, q_i$ $(1 \le i \le 3)$
satisfying the following equations:
$$
\begin{array}{ccccccc}
   n\ & = & n_1 & + & n_2 & + & n_3
 \\
   \rotatebox{90}{$=$} && \rotatebox{90}{$=$} 
                       && \rotatebox{90}{$=$} && \rotatebox{90}{$=$}
 \\
   p & = & p_1 & + & p_2 & + & p_3
 \\
   +   && +     && +     && +
 \\
   q & = & q_1 & + & q_2 & + & q_3
\end{array}
$$
We set
$$
J_i :=
  \begin{pmatrix}
   \;
     \begin{matrix}
            1 \\ &\ddots \\  &&1
     \end{matrix}
      \kern-3.2em\smash{\raisebox{4.8ex}{\rotatebox{-40}{$\overbrace{~~~~~~~~~~~~~}$}} }
             \kern-1.5em\smash{\raisebox{3.2ex}{$p_i$}}
  \\
    &\begin{matrix}
            0 \\ &\ddots \\  &&0
     \end{matrix}
      \kern-3.2em\smash{\raisebox{4.8ex}{\rotatebox{-40}{$\overbrace{~~~~~~~~~~~~~}$}} }
             \kern-1.5em\smash{\raisebox{3.2ex}{$q_i$}}
   \;
  \end{pmatrix}
\in  M(n_i, \mathbb{R} )
\quad (1 \le i \le 3).  
$$
Then 
$G_J \simeq U(p_1+p_2+p_3) \times U(q_1+q_2+q_3) = H$,
where $J$ is defined as in \eqref{eqn:J}.

Now, let us take a specific choice of
$X \in \mathfrak{u}(n)$ as follows:
for $z \in \mathbb{C}$,
\begin{equation*}
X_{12} := \begin{pmatrix}
             && 1 \\ & \mbox{\Large $0$} \\ 1
         \end{pmatrix}, \;
X_{23} := \begin{pmatrix}
             && 1 \\ & \mbox{\Large $0$} \\ 1
         \end{pmatrix}, \;
X_{13} := \begin{pmatrix}
             && 1 \\ & \mbox{\Large $0$} \\ z
         \end{pmatrix}.
\end{equation*}
This is possible because $n_i \ge 2$ $(i = 1,2,3)$.
Then the $(i,j)$-block of $Q := [X,J]$ amounts to
\begin{equation*}
Q_{12} 
 =  \begin{pmatrix}
                 && -1 \\  & \mbox{\Large $0$} \\ 1 
               \end{pmatrix}
, \;
Q_{23}
 =  \begin{pmatrix}
                 && -1 \\  & \mbox{\Large $0$} \\ 1 
               \end{pmatrix}
, \;
Q_{13} 
 = \begin{pmatrix}
                 && -1 \\  & \mbox{\Large $0$} \\ z
               \end{pmatrix}.
\end{equation*}
Associated to the loop $1 \to 2 \to 3 \to 1$,
we have
$$
A_{1231}(Q) = Q_{12} Q_{23} Q_{13}^*
   = \begin{pmatrix}
       && - \overline{z}  \\  & \mbox{\Large $0$} \\ 1
     \end{pmatrix}.
$$
Hence, 
$\det (\lambda I_{n_1} - A_{1231} (Q))
   = \lambda^{n_1} + \lambda^{n_1 -2} \overline{z}
   \notin \mathbb{R} [\lambda]
$
if we take $z \notin \mathbb{R}$.
Now,  Proposition~\ref{prop:notreal} follows from
Lemma~\ref{lem:machinenotreal}.
\end{proof}

\begin{proof}[Proof of Proposition \ref{prop:notreal2}]
With the notation as in \eqref{eqn:J},
we define $J$ by setting $k=4$ and
$$
J_1 = \begin{pmatrix}1\\ &\mbox{\Large $0$}\end{pmatrix}, \;
J_2 = \begin{pmatrix}1\\ &\mbox{\Large $0$}\end{pmatrix}, \;
J_3 = 0, \;
J_4 = 0.
$$
Then the isotropy subgroup $G_J$ is conjugate to
$H \simeq U(2) \times U(n-2)$
by an element of $G' = O(n)$.
We consider the loop
$1 \to 3 \to 2 \to 4 \to 1$.
Then it follows from \eqref{eqn:Pe}
 that 
\begin{align*}
   A_{13241} (Q) 
  ={}& Q_{13} Q_{23}^* Q_{24} Q_{14}^*
\\
  ={}& (-J_1 X_{13}) (-J_2 X_{23})^* (-J_2 X_{24}) (-J_1 X_{14})^*
\\
  ={}& 
     \begin{pmatrix}
                     (\vec a, \vec b) (\vec d, \vec c) & 0 & \cdots & 0
     \\
                      0                                &   &        &
     \\
                      \vdots                      &    & \mbox{\Large $0$} &
     \\
                      0                                 &   &        & 
     \end{pmatrix},
\end{align*}
where $\vec a, \vec b \in \mathbb {C}^{n_3}$
 denote the first row vectors of $X_{13}$, $X_{23}$
 and $\vec c, \vec d \in \mathbb {C}^{n_4}$
 denote the first row vectors of $X_{14}$, $X_{24}$,
 respectively.  
Since $n_1$, $n_2$, $n_3$ and $n_4$ are positive integers,
 we can find $X \in \mathfrak {u}(n)$
 such that $(\vec a, \vec b)(\vec d, \vec c) \not \in \mathbb {R}$.  

Then
$$
     \det (\lambda I_{n_1} 
      - A_{13241} (Q))
     = \lambda^{n_1 -1} (\lambda - (\vec{a},\vec{b})(\vec{d},\vec{c}))
                              \notin \mathbb {R}[\lambda].
$$
Thus Proposition \ref{prop:notreal2} follows from
Lemma~\ref{lem:machinenotreal}.  
\end{proof}

\begin{proof}[Proof of Proposition \ref{prop:notreal3}]

Without loss of generality,
 we may and do assume $1 \le n_1 \le n_2 \le \cdots \le n_k$.  
We give a proof according to the following three cases:

Case 1) \quad $k \ge 5$.

Case 2) \quad $k = 4$ and $n_3 > 1$.

Case 3) \quad $k = 4$ and $n_3 = 1$.

\medskip

Case 1) \enspace
Suppose $k \ge 5$.
Then,
 $n_5 + \cdots + n_k >1$.
In fact, if it were not the case,
we would have $k=5$ and $n_5=1$,
 which would imply $n_1 =n_2 = \cdots =n_5=1$ and $n=5$.
But, this contradicts to the assumption $n=p+q \ge 3+3=6$.  
Now, we apply Proposition \ref{prop:notreal}
 to $L'=U(n_1 + n_2) \times U(n_3+ n_4) \times U(n_5 + \cdots + n_k)$ ($\supset L$),
 and  conclude 
$$
G \supsetneqq L' G' H \supset L G' H.
$$

Case 2) \enspace
Suppose
 $k=4$ and  $n_3 >1$.
Then, $2 \le n_3 \le n_4$ and $2 \le n_1+n_2$.  
Then apply Proposition \ref{prop:notreal} to 
$L'=U(n_1+ n_2) \times U(n_3) \times U(n_4)$
 and we conclude that 
$$
G \supsetneqq L' G' H \supset L G' H.
$$

Case 3) \enspace
Suppose $k=4$ and $n_3 =1$.
Then,
\begin{equation}\label{eqn:kn}
     k=4 \quad\mbox{and}\quad n_1=n_2=n_3=1.  
\end{equation}
In the setting \eqref{eqn:J},
we define $J$ by setting
$$
   J_1=J_2=J_3:=(1) \in M(1,\mathbb {R})
   \quad\mbox{and}\quad
   J_4:=
  \begin{pmatrix}
   \;
     \begin{matrix}
            1 \\ &\ddots \\  &&1
     \end{matrix}
      \kern-3.2em\smash{\raisebox{4.8ex}
                {\rotatebox{-40}{$\overbrace{~~~~~~~~~~~~~}$}}}
             \kern-1.5em\smash{\raisebox{3.2ex}{{\scriptsize{$p-3$}}}}
  \\
    &\kern-1em{\begin{matrix} 0 \\ &\ddots \\  &&0 \end{matrix}}
      \kern-3.2em\smash{\raisebox{4.8ex}{\rotatebox{-40}{$\overbrace{~~~~~~~~~~~~~}$}} }
             \kern-1.5em\smash{\raisebox{3.2ex}{{\footnotesize{$q$}}\,\,}}
   \;
  \end{pmatrix}.  
$$
Then, $G_J \simeq U(p) \times U(q) = H$.

Given $Y = (v_1,\dots,v_p) \in M(q,p;\mathbb {C})$
$(v_1,\dots,v_p \in \mathbb{C}^q)$,
we define $X \in \mathfrak u(n)$
 by
$$
   X:=\begin{pmatrix} 0 & -Y^* 
      \\
                       Y & 0\end{pmatrix}.  
$$
Then, 
$Q:= [X,J] = \begin{pmatrix} 0 & Y^* \\ Y & 0 \end{pmatrix}$.
Associated to the partition 
$n=n_1+ n_2+n_3+n_4$ $(=1+1+1+(n-3))$
and the loop $1 \to 4 \to 2 \to 4 \to 3 \to 4 \to 1$,
\begin{align}
A_{1424341} (Q)
&= Q_{14} Q_{24}^* Q_{24} Q_{34}^* Q_{34} Q_{14}^*
\nonumber
\\
&= \begin{pmatrix}0 \\ v_1 \end{pmatrix}^*
   \begin{pmatrix}0 \\ v_2 \end{pmatrix}
   \begin{pmatrix}0 \\ v_2 \end{pmatrix}^*
   \begin{pmatrix}0 \\ v_3 \end{pmatrix}
   \begin{pmatrix}0 \\ v_3 \end{pmatrix}^*
   \begin{pmatrix}0 \\ v_1 \end{pmatrix}
\nonumber
\\
&= (v_2,v_1)(v_3,v_2)(v_1,v_3),
\label{eqn:V3}
\end{align}
where $0$ denotes the zero vector in $\mathbb{C}^{p-3}$.
Since $p \ge 3$ and $q \ge 3$ $(\ge 2)$,
we can take
$$
Y = (v_1,v_2,\dots,v_p) = \begin{pmatrix}
                            \begin{matrix} 1&1&1 \\ 1&z&0 \end{matrix}
                            & & \mbox{\Large $0$} \;
                          \\[3\medskipamount]
                             \mbox{\Large $0$}   & & \mbox{\Large $0$} \;
                          \end{pmatrix},
$$
so that $(v_2,v_1)(v_3,v_2)(v_1,v_3) = 1+z$.
Thus,
$$
\det (\lambda I_1 - A_{1424341} ([X,J]))
  = \lambda - 1 - z \notin \mathbb{R} [\lambda]
$$
if $z \notin \mathbb{R}$.
Hence Proposition follows from Lemma \ref{lem:machinenotreal}.
\end{proof}
Hence, the proof of Theorem~\ref{thm:A} is now completed.

\section{Visible actions on generalized flag varieties}
\label{sec:visible}

Suppose a Lie group $L$ acts holomorphically on a connected complex
manifold $D$ with complex structure $J$.

\begin{definition}
[{\rm{see {\cite[Definitions 2.3, 3.1.1 and 3.3.1]
{xkrims}}}}]
\label{def:visible}

The action is \textit{previsible} if there exists a totally real
submanifold $S$ such that
\begin{equation}\label{eqn:V1}
\text{$D' := L \cdot S$ \
      is open in $D$.}
\end{equation}
The previsible action is \textit{visible} if
\begin{equation}\label{eqn:V2}
J_x(T_xS) \subset T_x(L \cdot x)
\quad\text{for generic $x \in S$ ($J$-transversality)},
\end{equation}
and is {\it{strongly visible}}
 if there exists an anti-holomorphic
diffeomorphism
$\sigma$ of $D'$
such that
\begin{align}
&\sigma|_S = \operatorname{id},
\label{eqn:v3}
\\
&\text{$\sigma$ preserves each $L$-orbit on $D'$.}
\label{eqn:V4}
\end{align}
\end{definition}

A strongly visible action is visible (\cite[Theorem 14]{xkrims}).
Furthermore, we have:
\begin{lemma}
\label{lem:compati}
Suppose there exists an automorphism $\tilde{\sigma}$ of $L$ such that
\begin{equation}
\label{eqn:compati}
\sigma(g \cdot x) = \tilde{\sigma}(g) \cdot \sigma(x)
\quad (g \in L, \ x \in L).
\end{equation}
Then, a previsible action satisfying \eqref{eqn:v3} is strongly
visible.
\end{lemma}

\begin{proof}
Any $L$-orbit on $D'$
is of the form $L \cdot x$
 for some $x \in S$.  
Then, 
 by \eqref{eqn:v3} and \eqref{eqn:compati},
 we have
$
\sigma(L \cdot x) = \tilde{\sigma}(L) \cdot \sigma(x)
= L \cdot x.
$
Hence,
 the condition \eqref{eqn:V4} is fulfilled.  
\end{proof}

Now, let us consider the setting of Theorem~\ref{thm:A}.

\begin{example}\label{ex:compact sigma}
We set
$$
\tilde{\sigma}(g) := \overline{g}
\quad\text{for $g \in G = U(n)$.}
$$
Then $\tilde{\sigma}$ stabilizes subgroups $L$ and $H$ in the setting 
\eqref{eqn:1.1} in particular,
induces a Lie group automorphism, denoted by the same letter
$\tilde{\sigma}$,
 of $L$,
and an anti-holomorphic diffeomorphism,
denoted by $\sigma$, 
of the homogeneous space
$$
G/H \simeq \mathcal{B}_{m_1,\dots,m_l}(\mathbb{C}^n).  
$$
The $L$-action on $G/H$ satisfies
 the compatibility condition
 \eqref{eqn:compati}.
Now, 
 we consider the totally real submanifold
$
S := \mathcal{B}_{m_1,\dots,m_l}(\mathbb{R}^n)
$ 
 in $G/H$.  
Since $\tilde{\sigma} = \operatorname{id}$ on
$G' = O(n)$ and $S \simeq G'/G' \cap H$, we have
 $\sigma|_S = \operatorname{id}$.  
Therefore, 
 if $S$ meets every $L$-orbit on $G/H$,
then the $L$-action on $G/H$
 is not only previsible
 by definition 
 but also is strongly visible
 by Lemma \ref{lem:compati}.  
Hence, 
 the assertion in Remark \ref{rem:1.5.2}
 is proved.  
\end{example}

\section{Applications to representation theory}
\label{sec:concl}

This section gives
 a flavor of some applications
 of Theorem~\ref{thm:A}
 to multiplicity-free theorems
 in representation theory.  

In \cite{mfbdle} (see also \cite{mf-korea} and
\cite[Theorem~2]{xkrims}),
we proved that the multiplicity-free property propagates from fibers
to spaces of holomorphic sections of equivariant holomorphic
bundles under a certain geometric condition.
The key assumption there is strongly visible
actions
 (Definition \ref{def:visible}) on base spaces.

First of all, we observe that
 one dimensional representations are
 obviously irreducible, 
 and therefore is multiplicity-free.  
Then, by \cite{mfbdle}, 
this multiplicity-free property propagates
 to the multiplicity-free property of the representation on the space
 ${\mathcal{O}}(G/H, {\mathcal{L}}_{\lambda})$
of holomorphic sections
 as an $L$-module
 for any $G$-equivariant holomorphic line bundle
 ${\mathcal{L}}_{\lambda} \to G/H$
 if $(G,L,H)$ satisfies (ii)$'$ 
(or any of the equivalent conditions)
 in Theorem~\ref{thm:A}.  
Then, by a theorem of Vinberg--Kimelfeld\cite{xvk},
 this implies that
$
     G/H 
     \simeq 
    {\mathcal{B}}_{m_1, \cdots, m_l}({\mathbb{C}}^n) 
$
 is a spherical variety of 
$
     L_{\mathbb{C}}
     \simeq
     GL(n_1,{\mathbb{C}}) \times \cdots \times 
     GL(n_k,{\mathbb{C}})
$,
namely, a Borel subgroup of $L_{\mathbb{C}}$ has an open orbit on
$\mathcal{B}_{m_1,\dots,m_l}(\mathbb{C}^n)$.

More generally,
 applying the propagation theorem of multiplicity-free
 property to higher dimensional fibers,
 we get from 
 Theorem~\ref{thm:A} 
 a new geometric proof
 of a number of
 multiplicity-free results including:
\begin{itemize}
\item
(Tensor product)\enspace
The tensor product of two irreducible representations $\pi_\lambda$
and $\pi_\mu$ of $GL(n,\mathbb{C})$ is multiplicity-free 
if the highest weight
$\mu \in \mathbb{Z}^n$ is of the form, 
\begin{equation}
\label{eqn:abpq}
 (\underbrace{a,\dots,a}_p , \underbrace{b,\dots,b}_q )
\quad
(a \ge b).
\end{equation}
for some $a,b \in \mathbb{Z}$ $(a\ge b)$ and some $p,q~(p + q =n)$
and if one
of the following conditions is satisfied:
\begin{itemize}
\item[1)\phantom{$'$}]
$\min(p,q) = 1$.
\item[1)$'$]
$a-b = 1$.
\item[2)\phantom{$'$}]
$\min(p,q) = 2$ and $\lambda$ is of the form
\begin{equation}
\label{eqn:xyz1}
( \underbrace{x,\dots,x}_{n_1} ,
  \underbrace{y,\dots,y}_{n_2} ,
  \underbrace{z,\dots,z}_{n_3} )
\quad
(x \ge y \ge z)
\end{equation}
\item[2)$'$]
$a-b = 2$ and $\lambda$ is of the form \eqref{eqn:xyz1}.
\item[3)\phantom{$'$}]
$\lambda$ is of the form \eqref{eqn:xyz1} satisfying
$$
\min (x-y, y-z, n_1, n_2, n_3) = 1.
$$
\end{itemize}
Stembridge \cite{xstemgl} gave a proof of this fact
using a combinatorial method by a case-by-case argument.
He proved also that this exhausts all multiplicity-free cases.

\item
(Restriction: $GL_n \downarrow GL_p \times GL_q$) \enspace
The representation $\pi_\lambda$ of $GL(n,\mathbb{C})$ is
multiplicity-free when restricted to 
$GL(p,\mathbb{C}) \times GL(q,\mathbb{C})$ 
if one of the three conditions (1), (2), and (3) is satisfied.

\item
(Restriction: $GL_n \downarrow GL_{n_1} \times
  GL_{n_2} \times GL_{n_3}$) \enspace
The irreducible representation $\pi_\mu$ of $GL(n,\mathbb{C})$ having
  highest weight $\mu$ of the form \eqref{eqn:abpq}
is multiplicity-free when restricted to
$GL(n_1,\mathbb{C}) \times GL(n_2,\mathbb{C}) \times
 GL(n_3,\mathbb{C})$
if
$\min(a-b,p,q) \le 2$.
\end{itemize}

It is noteworthy
 that the above three multiplicity-free results
(``triunity'') are obtained from a single geometric result, Theorem~\ref{thm:A}
(see \cite[Theorems~3.3, 3.4 and 3.6]{xkleiden}).
In particular, we get the equivalence (ii) $\Leftrightarrow$ (iv) in
Remark \ref{rem:spherical} by considering the $(G \times
 G)$-equivariant line bundle
$\mathcal{L}_\lambda \boxtimes \mathcal{L}_\mu
 \to G/L \times G/H$.

Besides, Theorem~\ref{thm:A} gives
 a new geometric proof of yet more
 multiplicity-free theorems 
 (see \cite[Theorems~19, 20]{xkrims}): 

\begin{itemize}
\item
($GL_p$--$GL_q$ duality)\enspace
The symmetric algebra $S(M(p,q;\mathbb{C})) \simeq S(\mathbb{C}^{pq})$
is multiplicity-free as a representation of
$GL(p,\mathbb{C}) \times GL(q,\mathbb{C})$.
\item
(Kac \cite{xkac})\enspace
$S(\mathbb{C}^{pq})$ is still
 multiplicity-free as a representation of
$GL(p-1,\mathbb{C}) \times GL(q,\mathbb{C})$.
\item
(Panyushev \cite{panyu})\enspace
Let $\mathcal {N}$ be a nilpotent orbit of $GL(n,\mathbb {C})$
 corresponding to a partition $(2^p \, 1^{n-2p})$.  
Then the representation of $GL(n,\mathbb {C})$ on the space
 $\mathbb {C}[\mathcal {N}]$ of regular functions
 is multiplicity-free.   
\end{itemize}

All of the examples discussed so far are finite dimensional.
As we saw in \cite{mf-korea}, we can also expect
from strongly visible actions yet more 
multiplicity-free theorems for infinite dimensional representations
for both continuous and discrete spectra.
Applications of Theorem~\ref{thm:A}
 (and its non-compact version)
  to infinite dimensional representations
  will be discussed in a future paper.


\begin{thebibliography}{99}

\bibitem{xbigo}
L. Biliotti and A. Gori,
Coisotropic and polar actions on complex Grassmannians,
Trans. Amer. Math. Soc.
\textbf{357} (2004), 1731--1751.


\bibitem{xborelji}
A. Borel and L. Ji,
Compactification of symmetric and locally symmetric spaces,
69--137,
In: Lie Theory, Unitary Representations and Compactifications of 
Symmetric Spaces (eds. J.-P. Anker and B. {\O}rsted),
Progr. Math. \textbf{229} (2005).

\bibitem{xfjann}
M. Flensted-Jensen,
Discrete series for semisimple symmetric spaces,
Ann. Math. \textbf{111} (1980), 253--311.

\bibitem{xhoogen}
B. Hoogenboom, 
Intertwining functions on compact Lie groups. 
CWI Tract, 5. 
Stichting Mathematisch Centrum,
 Centrum voor Wiskunde en Informatica, Amsterdam, 1984.
  

\bibitem{xkac}
V. Kac, 
 Some remarks on nilpotent orbits,
 J. Algebra \textbf{64} (1980),
 190--213.  

\bibitem{xkinv}
T. Kobayashi,
Invariant measures on homogeneous manifolds of reductive type,
J. reine angew. Math.
\textbf{490} (1997),
37--53.

\bibitem{xkmfjp}
T. Kobayashi,
Multiplicity-free branching laws for unitary highest weight modules,
Proceedings of the Symposium on Representation Theory
held at Saga, Kyushu 1997 (ed. K. Mimachi) (1997),
9--17.


\bibitem{xkleiden}
T. Kobayashi,
Geometry of multiplicity-free representations of $GL(n)$,
visible actions on flag varieties, and triunity, 
Acta Appl. Math.,
\textbf{81} (2004),
129--146.

\bibitem{xkrims}
T. Kobayashi,
Multiplicity-free representations and visible actions on complex
manifolds, 
Publ. RIMS,
\textbf{41} (2005), 497--549
(a special issue of Publications of RIMS commemorating the fortieth
anniversary of the founding of the Research Institute for Mathematical
Sciences).  

\bibitem{mfbdle}
T. Kobayashi,
Propagation of multiplicity-free property for holomorphic vector bundles,
preprint

\bibitem{xkvisible}
T. Kobayashi,
Visible actions on symmetric spaces,
preprint

\bibitem{mf-korea}
T. Kobayashi,
Multiplicity-free theorems of the 
restrictions of unitary highest weight modules
        with respect to reductive symmetric pairs, 
to appear in Progr. Math., Birkh\"auser.

\bibitem{xlittelmann}
P. Littelmann, 
 On spherical double cones,
 J. Algebra, \textbf{166} (1994), 142--157.  

\bibitem{xmatsu}
T. Matsuki,
Double coset decomposition of algebraic groups
arising from two involutions. I.,
J. Algebra, \textbf{175} (1995), 865--925.

\bibitem{xmatsuki}
T. Matsuki,  
Double coset decompositions of reductive Lie groups
 arising from two involutions, 
J. Algebra, \textbf{197} (1997),  49--91.


\bibitem{xmatsucpt}
T. Matsuki,  
Classification of two involutions
 on compact semisimple Lie groups and root systems,
J. Lie Theory, \textbf{12} (2002),  41--68.

\bibitem{xmaos}
T. Oshima and T. Matsuki,
Orbits on affine symmetric spaces under the action of
the isotropy subgroups, J. Math. Soc. Japan, 
\textbf{32} (1980), 399--414.

\bibitem{panyu}
D. Panyushev,
Complexity and nilpotent orbits,
Manuscripta Math.,
\textbf{83} (1994), 223--237.


\bibitem{xstemgl}
J. R. Stembridge,  
Multiplicity-free products of Schur functions, 
Ann. Comb., \textbf{5} (2001), 113--121.


\bibitem{xvk}
       \'E. B. Vinberg and B. N. Kimelfeld,
          Homogeneous domains on flag manifolds and spherical
                subgroups of semisimple Lie groups,
        Funct. Anal. Appl., \textbf{12} (1978),
        168--174



\end{thebibliography}
\end{document}